\documentclass[12pt]{article}
\usepackage[bottom]{footmisc}
\usepackage[backref=page]{hyperref}
\usepackage{amssymb, amsmath, amscd, mathtools}
\usepackage{graphicx}
\usepackage[utf8]{inputenc}
\usepackage[T1]{fontenc}
\usepackage{epsfig}
\usepackage{color}
\usepackage{subfig, color}
\usepackage{resizegather}

\usepackage{algorithm}
\usepackage{algorithmic}

\usepackage{etoolbox}
\patchcmd{\thebibliography}{\section*{\refname}}{}{}{}

\usepackage[utf8]{inputenc}
\usepackage{soul}

\begin{document}
\title{GAMA: A Novel Algorithm for Non-Convex Integer Programs}
\author{Hedayat Alghassi\footnote{halghassi@cmu.edu},\,  Raouf Dridi\footnote{rdridi@andrew.cmu.edu},\, Sridhar Tayur\footnote{stayur@cmu.edu}\\
\small CMU Quantum Computing Group\\
\small Tepper School of Business, Carnegie Mellon University, Pittsburgh, PA 15213\\
}

\maketitle

\begin{abstract}
Inspired by the decomposition in the hybrid quantum-classical optimization algorithm we introduced in \cite{alghassi_graver_2019}, we propose here a new (fully classical) approach to solving certain non-convex integer programs using Graver bases. This method is well suited when (a) the constraint matrix $A$ has a special structure so that its Graver basis can be computed systematically, (b) several feasible solutions can also be constructed easily and (c) the objective function can be viewed as many convex functions quilted together. Classes of problems that satisfy these conditions include Cardinality Boolean Quadratic Problems (CBQP), Quadratic Semi-Assignment Problems (QSAP) and Quadratic Assignment Problems (QAP). Our Graver Augmented Multi-seed Algorithm (GAMA) utilizes augmentation along Graver basis elements (the improvement direction is obtained by comparing objective function values) from these multiple initial feasible solutions.   We compare our approach  with a best-in-class commercially available solver (Gurobi). Sensitivity analysis indicates that the rate at which GAMA slows down as the problem size increases is much lower than that of Gurobi. We find that for several instances of practical relevance, GAMA not only vastly outperforms in terms of time to find the optimal solution (by two or three orders of magnitude), but also finds optimal solutions within minutes when the commercial solver is not able to do so in 4 or 10 hours (depending on the problem class) in several cases. 
\end{abstract}

\textbf{Keywords}: Graver basis, test sets, non-linear non-convex integer programming, contingency tables, (0,1)- matrices, computational testing.
\newpage
\tableofcontents 
\newpage

\section{Introduction}
Many hard to solve practical non-linear integer programming problems have a specially structured linear constraint matrix. We study an important subset of these classes--including Cardinality Boolean Quadratic Problems (CBQP), Quadratic Semi-Assignment Problems (QSAP), and Quadratic Assignment Problems (QAP)--which have two features: (1) their Graver bases can be calculated systematically, and (2) multiple feasible solutions, which are uniformly spread out in the space of solutions, can be likewise systematically constructed. The hardness of such problems, then, stems from the non-convexity of their nonlinear cost function, and not from the Graver basis or the ability to find feasible solutions. In this paper, we target such problems and solve them using a novel approach, GAMA: Graver Augmented Multi-seeded Algorithm. 
\\~~

Our approach here is inspired by an earlier work (\cite{alghassi_graver_2019}) in which we introduced a  decomposition for a broader class (there was no restriction on $A$ but we required a convex objective function), with three components: 
\begin{itemize}
\item [(1)] Computing a partial or complete Graver basis of $A$ (which utilized outputs from a D-Wave quantum computer, a stand-in for any Ising solver);
\item [(2)] a set of feasible solutions (which also required a  call to D-Wave); and
\item [(3)] a parallel augmentation procedure, starting at each of the feasible solutions from (2), using the (partial or complete) Graver basis that was obtained in (1).
\end{itemize}
An essential concept that was exploited there, which we will continue to use here, was to separate the objective function from the constraints in the decomposition. It is known that given the Graver basis of  a problem’s integer matrix, any convex non-linear problem can be globally optimized with polynomial number of moves--augmentations--from any arbitrary initial feasible solution \cite{onn_nonlinear_2010}. What is novel in this current paper is the recognition that for a wide class of hard problems, the matrix $A$ has a special structure that allows us to obtain (a) Graver basis elements and (b) many feasible solutions that are spread out, by classical methods, and that are simple enough to be systematically algorithmized. Therefore, the steps (1) and (2) that previously needed calls to D-Wave no longer do.

~~\\
Furthermore, suppose the \textit{non-convex objective function} can be viewed as \textit{many convex\footnote{We have an expanded notion of convex; see Section \ref{subsec:costcategories}.} functions stitched together} (like a quilt). Thus, the entire feasible solution space can be seen as a collection of parallel subspaces, each with a convex objective function. If we have the Graver basis for the constraint matrix, and a feasible solution in every one of these sub-regions, putting this all together, an algorithm that can find the optimal solution is as follows:
\begin{itemize}
    \item Find the Graver basis.
    \item Find a  number of feasible solutions, spread out so that there is at least one feasible solution in each of the sub-regions that has a convex objective function.
    \item Augment along the Graver basis from each of the feasible solutions ("seeds") until you end up with a number of  local optimal solutions (one for each seed).
    \item Choose the best from among these local optimal solutions.
\end{itemize}
In this paper, we compare the speed of reaching optimal solutions with a best-in-class integer optimization solver, Gurobi\textsuperscript{\textregistered}\cite{gurobi}. We do not discuss the worst-case complexity of GAMA. Our goal is to solve instances from industry that are not solvable by commercially best-in-class solvers and understand (numerically) why GAMA works so well when it does.

\section{Recap: Graver Bases and their Use in Non-linear Integer Optimization}\label{sec:background}

This section is repeated verbatim from \cite{alghassi_graver_2019} for ease of access. Let ${f:{\mathbb{R}^n} \to \mathbb{R}}$  be a real-valued function.
We want to solve the general non-linear integer optimization problem:
\begin{equation}\label{eq:gennonlin}
{(IP)_{A,b,l,u,f}} : \left\{ {\begin{array}{*{20}{l}}
  {\begin{array}{*{20}{c}}
  {\min }&{f(x)},&& 
\end{array}} \\ 
  {\begin{array}{*{20}{c}}
  {Ax = b},&&{l \leqslant x \leqslant u},&&{x,l,u \in {\mathbb{Z}^n}} 
\end{array}} \\ 
  {\begin{array}{*{20}{c}}
  {A \in {\mathbb{Z}^{m \times n}}},&&{b \in {\mathbb{Z}^m}} 
\end{array}} 
\end{array}} \right.
\end{equation}
One approach to solving such problem is to use an augmentation procedure: start from an initial feasible solution (which itself can be hard to find) and take an improvement step ({\it augmentation}) until one reaches the optimal solution. Augmentation procedures such as these need {\it test sets} or {\it optimality certificates}:  so it either declares the optimality of the current feasible solution or provides direction(s) towards better solution(s). Note that it does not matter from which feasible solution one begins, nor the sequence of improving steps taken: the final stop is an optimal solution.
\paragraph{Definition 1.} 
A set  $\mathcal{S} \in {\mathbb{Z}^n}$ is called a test set or optimality certificate if for every non-optimal but feasible solution ${x_0}$ there exists  $t \in \mathcal{S}$ and  $\lambda  \in {\mathbb{Z}_ + }$ such that  ${x_0} + \lambda t$ is feasible and $f\left( {{x_0} + \lambda t} \right) < f\left( {{x_0}} \right)$. The vector $t$ (or $\lambda t$)  is called the  \textit{improving} or \textit{augmenting} direction.
\\~~
If the optimality certificate is given, any initial feasible solution ${x_0}$ can be augmented to the optimal solution.  If $\mathcal{S}$ is finite, one can enumerate over all $t \in \mathcal{S}$ and check if it is augmenting (improving). If $\mathcal{S}$ is not practically finite, or if all elements $t \in \mathcal{S}$ are not available in advance, it is still practically enough to find a subset of $\mathcal{S}$ that is feasible and augmenting.  
\begin{figure}[H]
\centering
\includegraphics[width=3.5cm]{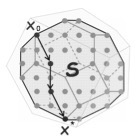}
\caption{Augmenting from initial to optimal solution \cite{onn_nonlinear_2010}.}
\label{fig:augmentation}
\end{figure}

~~\\
In the next section, we discuss the Graver basis of an integer matrix   $A \in {\mathbb{Z}^{m \times n}}$ which is known to be an  {optimality certificate}. 

\subsection{Mathematics of Graver Bases}
First, on the set ${\mathbb{R}^n}$,  we define
the following partial order:

\paragraph{Definition 2.} 
Let  $x,y \in {\mathbb{R}^n}$.
We say $x$  is {\it conformal} to $y$, 
written  $x \sqsubseteq y$, when   ${x_i}{y_i} \geqslant 0$ ($x$ and $y$ lie on the same orthant) and $\left| {{x_i}} \right| \leqslant \left| {{y_i}} \right|$  for $i = 1,...,n$. 
Additionally, a sum $u = \sum\limits_i {{v_i}}$ is called {\it conformal}, if ${v_i} \sqsubseteq u$ for all $i$ ($u$ \textit{majorizes} all $v_i$).  


~~\\
Suppose $A$ is a matrix in $ {\mathbb{Z}^{m \times n}}$. Define:   
\begin{equation}\label{eq:kernel}
{\mathcal{L}^*}(A) = \left\{ {x\left| {\begin{array}{*{20}{c}}
  {Ax = {\mathbf{0}},}&{x \in {\mathbb{Z}^n}\begin{array}{*{20}{c}}
  ,&{A \in {\mathbb{Z}^{m \times n}}} 
\end{array}} 
\end{array}} \right.} \right\}\backslash \left\{ {\mathbf{0}} \right\}.
\end{equation}
The notion of the Graver basis was first introduced in \cite{graver_foundations_1975} for integer linear programs (ILP):
\paragraph{Definition 3.} 
The Graver basis  of   integer matrix $A$  is defined to be the finite set of $\sqsubseteq $ minimal elements (\textit{indecomposable} elements) in the lattice  ${\mathcal{L}^*}(A)$. We denote 
by $\mathcal{G}(A) \subset {\mathbb{Z}^n}$ the Graver basis of $A$. 

~~\\
The following proposition summarizes the properties of Graver bases that are relevant to our setting. 
\paragraph{Proposition 1.} 
The following statements are true: 
\begin{itemize}
\item[(i)] Every vector $x$  in the lattice   ${\mathcal{L}^*}(A)$   is a conformal sum of the Graver basis elements.
\item[(ii)]  Every vector $x$  in the lattice   ${\mathcal{L}^*}(A)$   can be written as 
$
{x = \sum\limits_{i = 1}^t {{\lambda _i}} {g_i}}$ for some $ \lambda _i \in {\mathbb{Z}_ + }$ and $ 
   g_i \in \mathcal{G}(A) $. The upper bound on the number of Graver basis elements required ($t$) (called \textit{integer Caratheodory number}) is $\left( {2n - 2} \right)$.
   
\item[(iii)] A Graver basis is a test set (optimality certificate) for ILP and several nonlinear convex forms (see \ref{subsec:costcategories}). That is, a point ${x^*}$ is optimal for the optimization problem ${(IP)_{A,b,l,u,f}}$, if and only if there are no ${g_i} \in \mathcal{G}(A)$ such that ${x^*} + {g_i}$ is better than  ${x^*}$.

\item[(iv)] For any $g \in \mathcal{G}(A)$, an upper bound on the norm of Graver basis elements is given by
\begin{equation}\label{eq:graverupper}
{\left\| g \right\|_\infty } \leqslant \left( {n - r} \right)\Delta \left( A \right) \mbox{ and } \  {\left\| g \right\|_1 } \leqslant \left({n-r} \right)\left({r+1} \right)\Delta \left( A \right)
\end{equation}
where $r = rank(A) \leqslant m$  and $\Delta \left( A \right)$ is the maximum absolute value of the determinant of a square submatrix of $A$. 
\end{itemize}

 
 
\subsection{Applicability of Graver Bases as Optimality Certificates}\label{subsec:costcategories}
Beyond integer linear programs (ILP) with a fixed integer matrix, Graver bases as optimality certificates have now been generalized to include several nonlinear objective functions: 
\begin{itemize}
    \item Separable convex minimization \cite{murota_optimality_2004}: $\begin{array}{*{20}{c}}
  {\min }&{\sum\nolimits_i {{f_i}({c_i^Tx})} } 
\end{array}$ with $f_i$  convex.

    \item Convex integer maximization (weighted) \cite{de_loera_convex_2009}: $\begin{array}{*{20}{c}}
  {\max }&{f({\mathbf{W}}x)}&{,{\mathbf{W}} \in {\mathbb{Z}^{d \times n}}} 
\end{array}$ with $f$ convex on $\mathbb{Z}^d$.
    
    \item  Norm $p$ (nearest to $x_0$) minimization  \cite{hemmecke_polynomial_2011}: $\begin{array}{*{20}{c}}
  {\min }&{{{\left\| {x - {x_0}} \right\|}_p}} 
\end{array}$. 
    
    \item  Quadratic minimization \cite{murota_optimality_2004, lee_quadratic_2012}: $\begin{array}{*{20}{c}}
  {\min }&{{x^T}Vx} 
\end{array}$ where $V$ lies in the dual of quadratic Graver cone of $A$ 
    
    \item Polynomial minimization \cite{lee_quadratic_2012}: $\begin{array}{*{20}{c}}
  {\min }&{P(x)} \end{array}$ where  $P$ is a polynomial of degree $d$, that lies on cone $K_d(A)$, the dual of $d^{th}$ degree Graver cone of $A$.
 
\end{itemize}
 It has been shown that only polynomially many augmentation steps are needed to solve such  minimization problems \cite{hemmecke_polynomial_2011} \cite{de_loera_convex_2009}. 
 \\~~
 
 In the rest of this paper, we loosely call all of the above mentioned cost categories as \textit{convex}.
 Graver did not provide an algorithm for computing Graver bases; Pottier \cite{pottier_euclidean_1996} and Sturmfels \cite{sturmfels_variation_1997} provided such algorithms.  The use of test sets to study integer programs is a vibrant area of research. Another collection of test sets is based on Groebner bases. See \cite{conti_buchberger_1991,  tayur_algebraic_1995,bigatti_computing_1999, hosten_grin:_1995} and \cite{bertsimas_new_2000} for several different examples of the use of Groebner bases for integer programs.

\section{Graver Bases for Some Structured Matrices}\label{sec:structured}
Here we construct the Graver basis of four categories of structured matrices that occur frequently as subsets of constraints in nonlinear integer programming problems. 
\subsection{Graver basis of \texorpdfstring{$\mathbf{1}_{k}^T$}{Lg}}
The matrix $A = \mathbf{1}_{k}^T$ (a row vector where all $k$ elements are 1) is unimodular\footnote{The determinant of every square submatrix of a unimodular matrix is $0$ or $\pm 1$.}. Therefore, the Graver basis elements contain only $\{-1, 0, +1\}$ values.
~~\\

For any ${g} \in \mathcal{G}(A)$, the number of nonzero elements of $g$ should be even with equal numbers of $+1$ and $-1$.
Starting from vectors with $2$ nonzero elements, we can have $\left( {\begin{array}{*{20}{c}}
  k \\ 
  2 
\end{array}} \right)$ vectors of the form $g = {e_i} - {e_j}$, for all $i=1,2,..., k$ and $j=i+1, ..., k$. We cannot have Graver elements with $4$ or more elements, since each of them is the positive sum of the $2$ nonzero vectors (also each of $e_i - e_j$ elements lying on a $(k-1)$ dimensional hyperplane is $\sqsubseteq -$minimal to any $e_i - e_j + e_k - e_l$ vector lying on $(k-3)$ dimensional hyperplane).
Therefore, all of $\pm (e_i - e_j)$ vectors construct the Graver basis of $A = \mathbf{1}_{k}^T$.
\begin{equation}\label{eq:graverone}
\begin{array}{*{20}{c}}
  {\mathcal{G}({{\mathbf{1}}_k^T}) = \left\{ { \pm ({e_i} - {e_j})} \right\}}&{\left\{ {\begin{array}{*{20}{c}}
  {i = 1, \ldots ,k} \\ 
  {j = i+1, \ldots ,k} 
\end{array}} \right.} 
\end{array}
\end{equation}

\subsection{Graver basis of \texorpdfstring{${I_n \otimes {\mathbf{1}_k}^T}$}{Lg}}
For any block diagonal matrix of the form,
${\mathbf{A}} = {I_n} \otimes A = \left[ {\begin{array}{*{20}{c}}
  A&{\mathbf{0}}& \cdots &{\mathbf{0}}&{\mathbf{0}} \\ 
  {\mathbf{0}}&A& \cdots &{\mathbf{0}}&{\mathbf{0}} \\ 
   \vdots & \vdots & \ddots & \vdots & \vdots  \\ 
  {\mathbf{0}}&{\mathbf{0}}& \cdots &A&{\mathbf{0}} \\ 
  {\mathbf{0}}&{\mathbf{0}}& \cdots &{\mathbf{0}}&A 
\end{array}} \right],$
 one can easily observe that
\[\mathcal{G}({\mathbf{A}}) = {I_n} \otimes \mathcal{G}(A) = \left[ {\begin{array}{*{20}{c}}
  {\mathcal{G}(A)}&{\mathbf{0}}& \cdots &{\mathbf{0}}&{\mathbf{0}} \\ 
  {\mathbf{0}}&{\mathcal{G}(A)}& \cdots &{\mathbf{0}}&{\mathbf{0}} \\ 
   \vdots & \vdots & \ddots & \vdots & \vdots  \\ 
  {\mathbf{0}}&{\mathbf{0}}& \cdots &{\mathcal{G}(A)}&{\mathbf{0}} \\ 
  {\mathbf{0}}&{\mathbf{0}}& \cdots &{\mathbf{0}}&{\mathcal{G}(A)}
\end{array}} \right].\]
~~\\
Therefore,
\begin{equation} \label{eq:gravermanyone}
\mathcal{G}({I_n} \otimes {\mathbf{1}}_k^T) = {I_n} \otimes \mathcal{G}({\mathbf{1}}_k^T)
\end{equation}
and $\mathcal{G}({\mathbf{1}}_k^T)$ can be acquired from equation (\ref{eq:graverone}).
\subsection{Graver basis of \texorpdfstring{${{\mathbf{1}_n}^T} \otimes I_k$}{Lg}}
In matrices of the form 
${\mathbf{A}} = {\mathbf{1}}_n^T \otimes {I_k} = \overbrace {\left[ {\begin{array}{*{20}{c}}
  {{I_k}}&{{I_k}}& \cdots &{{I_k}} 
\end{array}} \right]}^{\begin{array}{*{20}{c}}
  n&{matrices} 
\end{array}},$
there are $n$ number of $1$'s in each row; each consecutive pair of $1$'s is spaced by $k-1$ number of $0$'s, and there is only one $1$ in each column. The Graver basis of each row, then, is the Graver basis of $n$ multiple of $1$'s, but each adjacent pair of elements in it should be spaced apart by $k$ elements. This has to repeat for all of the $k$ rows. We can represent the $k$ element spacing for all $k$ rows by $\otimes {I_k}$. Therefore,
\begin{equation} \label{eq:gravermanyeye}
\mathcal{G}({\mathbf{A}}) = \mathcal{G}({\mathbf{1}}_n^T) \otimes {I_k}
\end{equation}
and $\mathcal{G}({\mathbf{1}}_n^T)$ is acquired from equation (\ref{eq:graverone}).

\subsection{Graver basis of combination: \texorpdfstring{$\left( {{\mathbf{1}}_n^T \otimes {I_k}} \right) \oplus \left( {{I_n} \otimes {\mathbf{1}}_k^T} \right)$}{Lg}}
This structured combination is  called the \textit{generalized Lawrence configuration}.
\subsubsection{Generalized Lawrence configuration and \texorpdfstring{$n$}{Lg}-fold product matrices}
 The Lawrence lifting of an $m \times k$ matrix $A$ is the enlarged matrix 
$\Lambda (A) = \left[ {\begin{array}{*{20}{c}}
  {{I_k}}&{{I_k}} \\ 
  A&{\mathbf{0}} 
\end{array}} \right] \in {\mathbb{Z}^{(k + m) \times 2k}}$ where $\mathbf{0}$ is the $m \times k$ matrix of all zeros and $I_k$ is the $k \times k$ identity matrix. The generalization of the Lawrence lifting is of the form
\[{\Lambda _n}(A) = \left[ {\begin{array}{*{20}{c}}
  {{I_k}}&{{I_k}}& \cdots &{{I_k}}&{{I_k}} \\ 
  A&{\mathbf{0}}& \cdots &{\mathbf{0}}&{\mathbf{0}} \\ 
   \vdots &A& \ddots & \vdots & \vdots  \\ 
  {\mathbf{0}}& \vdots & \ddots &{\mathbf{0}}&{\mathbf{0}} \\ 
  {\mathbf{0}}&{\mathbf{0}}& \cdots &A&{\mathbf{0}} 
\end{array}} \right].\]
 ~~\\
 Graver bases  of ${\Lambda _n}(A)$ are finite (\cite{santos_higher_2003}, \cite{hosten_finiteness_2007}). A slightly modified version of this called $n$-fold matrices \cite{de_loera_n-fold_2008} appears frequently in integer programming. The $n$-fold matrix appears in many applications including high dimensional transportation problems and packing problems. Our $n$-fold matrix is constructed of \textit{ordered pair} $I_{k}, \mathbf{1}_k^T$ (${\mathbf{A}} = [I_k,\mathbf{1}_k^T]^n$). Given the Graver basis of fixed sized matrix $A$, there is a polynomial time algorithm that computes the Graver basis of the $n$-fold product matrix \cite{onn_nonlinear_2010}. 
\subsubsection{Systematic generation of Graver basis}\label{ss:combgraver}
If we had only the second term $\left( {{I_n} \otimes {\mathbf{1}}_k^T} \right)$, then based on equation (\ref{eq:gravermanyone}) the overall Graver basis would be the vector placement of each of $k$ sized vectors of $\mathcal{G}\left( {{\mathbf{1}}_k^T}\right)$ into $n$ bricks. The addition of $\left( {{\mathbf{1}}_n^T \otimes {I_k}} \right)$ block constrains the brick placement of those vectors, such that the sum of the placements should become zero. Keeping the vector placement aspect aside, we need to know which minimal and positive combination of elements of $\mathcal{G}\left( {{\mathbf{1}}_k^T}\right)$ reduces to zero. This is equivalent to the Graver basis of $\mathcal{G}\left( {{\mathbf{1}}_k^T}\right)$ in positive orthant, which is the Hilbert basis of $\mathcal{G}\left( {{\mathbf{1}}_k^T}\right)$. 
Therefore, we need to find the Hilbert basis of the Graver basis of $\mathbf{1}_k^T$ and do the liftings (brick placements) based on it.
~~\\

The following steps generate the Graver basis for $\left( {{\mathbf{1}}_n^T \otimes {I_k}} \right) \oplus \left( {{I_n} \otimes {\mathbf{1}}_k^T} \right)$:
\begin{itemize}
\item[(1)] Calculate the Hilbert basis of the Graver basis of $A = {\mathbf{1}}_k^T$ (see \ref{ss:Hilbert}):
\begin{equation}
\mathbf{H} = \mathcal{H}\left( {\mathcal{G}\left( {{\mathbf{1}}_k^T} \right)} \right).
\end{equation}

\item[(2)] Generate all $t$ to $n$ ($t=2,...,k$) liftings of $\mathcal{G}\left( {{\mathbf{1}}_k^T}\right)$ elements (bricks) based on Hilbert basis elements $\mathbf{H}$  and their variants (see \ref{ss:Liftings}).
\end{itemize}

\subsubsection{Hilbert basis of \texorpdfstring{$\mathcal{G}\left( {{\mathbf{1}}_k^T}\right)$}{Lg}}\label{ss:Hilbert}
One way to generate the Hilbert basis of $\mathcal{G}\left( {{\mathbf{1}}_k^T}\right)$ is using a completion procedure such as the Pottier algorithm \cite{pottier_euclidean_1996} and limiting it to the positive orthant. However, this method does not exploit the structure of $\mathcal{G}\left( {{\mathbf{1}}_k^T}\right)$ elements.

As evident in equation (\ref{eq:graverone}), the elements of $\mathcal{G}\left( {{\mathbf{1}}_k^T}\right)$ consists of $2\left( {\begin{array}{*{20}{c}}
  k \\ 
  2 
\end{array}} \right)$ vectors of size $k$, each containing one $+1$ and one $-1$ and $(k-2)$ zeros ($\{\pm ({e_i} - {e_j})\}$). We can model the matrix $\mathcal{G}\left( {{\mathbf{1}}_k^T}\right)$ as the \textit{incidence matrix} of a \textit{bidirectional complete graph} having $k$ nodes and $k(k-1)$ directional edges (two back and forth directional edges between each node pair). It can be shown that the set of all \textit{basic} (non-overlapping) directional cycles in this complete graph represents the Hilbert basis of $\mathcal{G}\left( {{\mathbf{1}}_k^T}\right)$.
~~\\

Therefore, we have $(2-1)!\left( {\begin{array}{*{20}{c}}
  k \\ 
  2 
\end{array}} \right)$ two-cycles, $(3-1)!\left( {\begin{array}{*{20}{c}}
  k \\ 
  3 
\end{array}} \right)$ three-cycles (directional triangles), $(4-1)!\left( {\begin{array}{*{20}{c}}
  k \\ 
  4 
\end{array}} \right)$ four-cycles, ..., and finally $(k-1)!\left( {\begin{array}{*{20}{c}}
  k \\ 
  k 
\end{array}} \right)$ $k$-cycles. The sum of these basic directional cycles, which is the cardinality of Hilbert basis results in:
\begin{equation}
Card\left( {\mathcal{H}\left( {\mathcal{G}({\mathbf{1}}_k^T)} \right)} \right) = 1!\left( {\begin{array}{*{20}{c}}
  k \\ 
  2 
\end{array}} \right) + 2!\left( {\begin{array}{*{20}{c}}
  k \\ 
  3 
\end{array}} \right) + 3!\left( {\begin{array}{*{20}{c}}
  k \\ 
  4 
\end{array}} \right) +  \cdots  + (k - 1)!\left( {\begin{array}{*{20}{c}}
  k \\ 
  k 
\end{array}} \right).
\end{equation}
One can thus construct each of the $t$-cycle indices of a $k$ complete graph and add them up for $t=2,...,k$.
~~\\
\subsubsection{Liftings and final Graver basis enumeration}\label{ss:Liftings}
Each of the Hilbert basis elements with $t$ nonzero terms (corresponding to a directional $t$-cycle), is a labeled construct, which creates $t!$ different possibilities. For each such possibility, we have $\left( {\begin{array}{*{20}{c}}
  n \\ 
  t 
\end{array}} \right)$ different lifting choices (combinations), which creates Graver elements and their negatives (by symmetry). Therefore, the total enumeration of the size $t$ Graver basis will be $\frac{{t!}}{2}\left( {\begin{array}{*{20}{c}}
  n \\ 
  t 
\end{array}} \right)$, excluding symmetries. 
~~\\
For each of the size $2$ Hilbert basis elements $\frac{2!}{2}\left( {\begin{array}{*{20}{c}}
  n \\ 
  2 
\end{array}} \right)$ different $2$-brick liftings, for each of the size $3$ elements $\frac{3!}{2}\left( {\begin{array}{*{20}{c}}
  n \\ 
  3 
\end{array}} \right)$ different $3$-brick liftings, for each of the size $4$ elements  $\frac{4!}{2}\left( {\begin{array}{*{20}{c}}
  n \\ 
  4 
\end{array}} \right)$ different $4$-brick liftings, ..., and finally for each of the size $k$ elements $\frac{k!}{2}\left( {\begin{array}{*{20}{c}}
  n \\ 
  k 
\end{array}} \right)$ different $k$-brick liftings. Combining the number of Hilbert basis elements in each size with all possible liftings, the total cardinality of the lifted Graver basis for matrix  $\left( {{\mathbf{1}}_n^T \otimes {I_k}} \right) \oplus \left( {{I_n} \otimes {\mathbf{1}}_k^T} \right)$ becomes:
\begin{equation}\nonumber
Card\left( {\mathcal{G}({\mathbf{A}})} \right) = 1!\frac{{2!}}{2}\left( {\begin{array}{*{20}{c}}
  k \\ 
  2 
\end{array}} \right)\left( {\begin{array}{*{20}{c}}
  n \\ 
  2 
\end{array}} \right) + 2!\frac{{3!}}{2}\left( {\begin{array}{*{20}{c}}
  k \\ 
  3 
\end{array}} \right)\left( {\begin{array}{*{20}{c}}
  n \\ 
  3 
\end{array}} \right) + 3!\frac{{4!}}{2}\left( {\begin{array}{*{20}{c}}
  k \\ 
  4 
\end{array}} \right)\left( {\begin{array}{*{20}{c}}
  n \\ 
  4 
\end{array}} \right) +  \cdots  + (k - 1)!\frac{{k!}}{2}\left( {\begin{array}{*{20}{c}}
  k \\ 
  k 
\end{array}} \right)\left( {\begin{array}{*{20}{c}}
  n \\ 
  k 
\end{array}} \right).
\end{equation}
Compressing,
\begin{equation}
Card\left( {\mathcal{G}({\mathbf{A}})} \right) = \frac{1}{2}\sum\limits_{t = 2}^k {\frac{{P(k,t)}}{t}} P(n,t)
\end{equation}
where $P$ is the permutation operator.
~~\\

\section{Computational Results}
\subsection{Algorithm}

{\centering
\begin{minipage}{1.0\linewidth}
\begin{algorithm}[H]
\caption{Graver Augmented Multi-seeded Algorithm (GAMA)}
\label{alg:NCGO}
\begin{algorithmic}[1]
\STATE {\bf inputs:} Matrix $A$, vector $b$, cost function $f(x)$, bounds $[l,u]$, as in Equation (\ref{eq:gennonlin})
\STATE {\bf output:} Global solution(s) $x^{*}$.
\STATE {\bf initialize:} terminated solutions set, $termSols = \{ \emptyset \}$.
\STATE Using any Graver extraction Algorithm input: $A$, extract $\mathcal G(A)$ (see section \ref{sec:structured})
\STATE Find multiple feasible solutions, satisfying ${l \leqslant x \leqslant u}$ (subsec. \textit{Feasible Solutions})
\FOR {any feasible solution $x = x_{0}$}
\WHILE{$g \in \mathcal G(A)$}
\IF{$l \leqslant(x+g)\leqslant u$ and $f(x+g) < f(x)$}
\STATE $x=x+g$
\ENDIF
\ENDWHILE
\STATE $termSols \leftarrow \left( {x,f(x)} \right)$
\ENDFOR
\RETURN ${x^*} = \left\{ {x\left| {f(x) = {{\min }_f}(termSols)} \right.} \right\}$
\end{algorithmic}
\end{algorithm}
\end{minipage}
}
\subsection{Cardinality Boolean Quadratic Programs (CBQP)}
 Cardinality Boolean Quadratic Programming (CBQP)  \cite{lima_solution_2017} is of the form:
\begin{equation}
\begin{array}{*{20}{c}}
  {(P1)}&{\min \left\{ {{\mathbf{c}^T}{\mathbf{x}} + {{\mathbf{x}}^T}Q{\mathbf{x}}:{\mathbf{1}}_n^T{\mathbf{x}} = b} \right\}} 
\end{array}
\end{equation}
where, 
$\begin{array}{*{20}{c}}
  {{{\mathbf{x}}^T} = \left[ {\begin{array}{*{20}{c}}
  {{x_1}}&{{x_2}}& \ldots &{{x_i}}& \ldots &{{x_n}} 
\end{array}} \right],}&&{{x_i} \in \mathbb{Z},} 
\end{array}$
 $\mathbf{c} \in {\mathbb{R}^n}$, and $Q \in {\mathbb{R}^{n \times n}}$.  Note that $Q$ is not necessarily a positive semidefinite matrix. Applications of this problem include edge-weighted graph problems \cite{billionnet_different_2005} as well as facility location problems \cite{bruglieri_annotated_2006}.
\subsubsection{Feasible solutions}
The Graver basis of $A = {\mathbf{1}}_n^T$ can be generated using equation (\ref{eq:graverone}). To solve the problem using our multi seeded approach, we need to generate many (say $l$) uniformly distributed feasible solutions of the linear constraint ${{\mathbf{1}}_n^Tx = b}$, $b \leqslant n$. 
The total number of solutions is $\left( {\begin{array}{*{20}{c}}
  n \\ 
  b 
\end{array}} \right)$, which can be large depending on the value of $b$. We generate only $l$ uniformly distributed number of them as initial feasible solutions. Each of the $l$ solutions is a size $n$ vector with $1$ placed in  $b$ locations with indices \textit{sampled uniformly at random} from $1$ to $n$ (without replacement), and $0$ placed in other locations.
\subsubsection{CBQP results}

All problems are also solved with the MIP solver of Gurobi\textsuperscript{\textregistered} Optimizer (latest version, 8.1), for comparison purposes\footnote{Installed on MATLAB R2014b using MacBookPro15,1: 6 Core 2.6 GHz Intel Core i7 processor, 32 GB 2400 MHz DDR4 RAM}. 
~~\\

We have used CBQP instances\footnote{https://sites.google.com/site/cbqppaper/} generated by \cite{lima_solution_2017}. There are $30$ CBQP instances of size $n=50$.
Setting the number of initial feasible points to be $l = 50$, our algorithm calculated the optimal solution for all 30 problems with an average time of $t_{av} = 1.052sec.$ In almost all cases, our method takes about about $1sec$. This is attributable to the fact that each augmentation path in all cases has a similar number of augmentation steps. With $\left( {\begin{array}{*{20}{c}}
  {50} \\ 
  2 
\end{array}} \right) = 1225$ Graver basis elements, each augmentation path takes about $0.02sec$ to terminate.

~~\\

\begin{figure}[H]
\centering
\includegraphics[width=17cm]{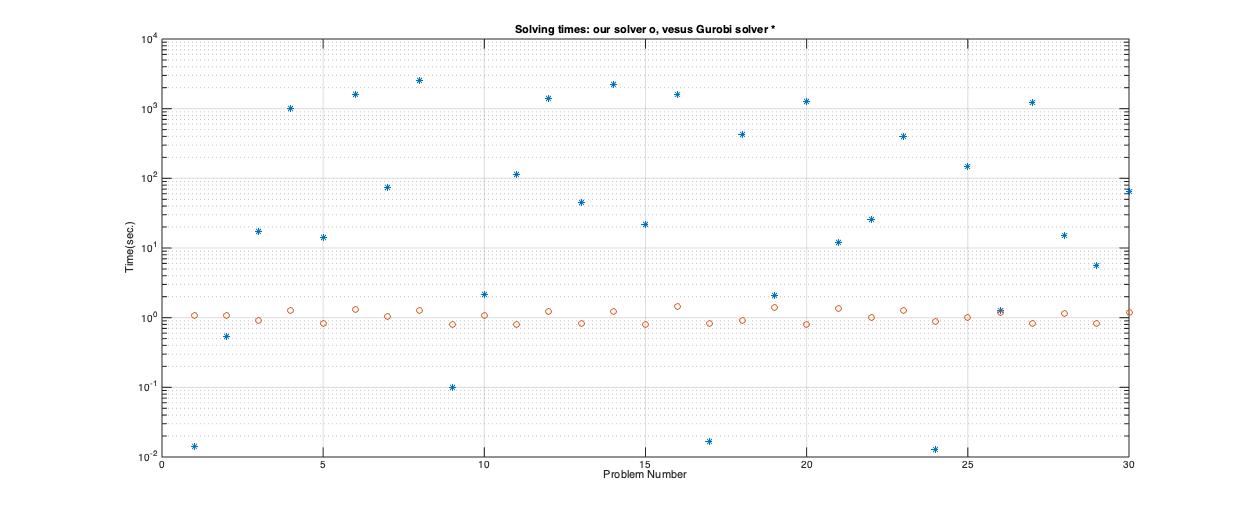}
\caption{Solving time of the 30 problem instances, GAMA 'o' vs. Gurobi '*'}
\label{fig:timesCBPQ}
\end{figure}

Using the Gurobi solver, in $5$ cases the problem was solved very quickly (in less than one second). At the other extreme, in $12$ cases the problems were very hard and took between $100sec$ to $2525sec$.  In the remaining 13 cases, the problems were solved between $1s$ and $100s$.The results are shown in Figure \ref{fig:timesCBPQ}.
~~\\

To summarize: our multi-seeded Graver based approach (GAMA) solved all $30$ instances (containing convex, easy nonconvex, and hard nonconvex instances) in a total of $30.55$ seconds; the same 30 instances took $~3.96\ hours$ for Gurobi. This is an average speedup of a factor of over $450$ for all samples. If we compare the performance only on the $12$ hard non-convex problems, we have a speedup of over $1000$.

\begin{figure}[H]
\centering
\includegraphics[width=10cm]{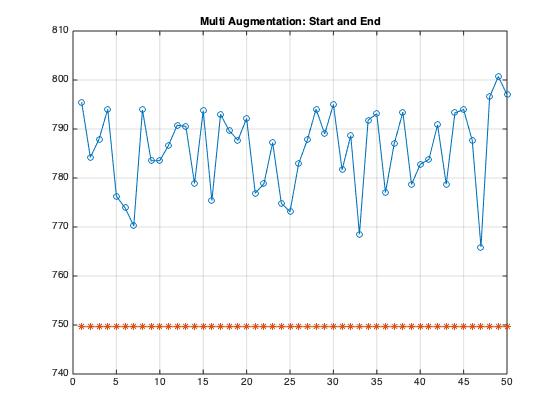}
\caption{All initial feasible solutions result in one global solution (Convex objective function).}
\label{fig:convexCBPQ}
\end{figure}

\begin{figure}[H]
\centering
\includegraphics[width=10cm]{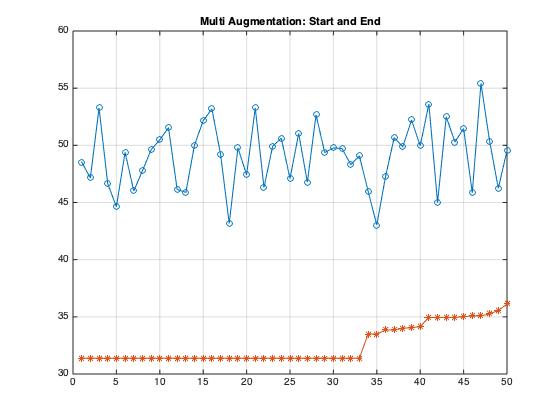}
\caption{Initial feasible solutions converge into a limited number of  terminal values (Non-Convex, but "easier").}
\label{fig:nonConvexLess}
\end{figure}

\begin{figure}[H]
\centering
\includegraphics[width=10cm]{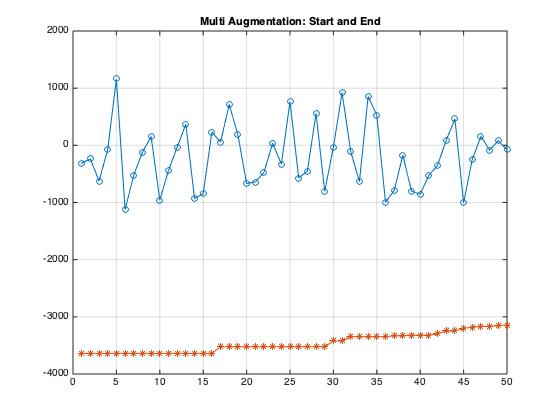}
\caption{Initial feasible solutions converge into a higher number of  different terminal values (Non-Convex,  "harder").}
\label{fig:nonConvexMore}
\end{figure}

~~\\
We want to understand when (and why) our approach does so much better when it does. We observed that the problem instances can be roughly categorized into three groups based on where all the augmentations ended up, beginning from the different starting feasible solutions.
\begin{itemize}
    \item All initial feasible solutions terminated at one global solution  (Figure \ref{fig:convexCBPQ}). This is the case when the objective function is convex\footnote{Recall that we call the  cost categories discussed in \ref{subsec:costcategories} \textit{convex}.}.
    \item The number of different terminal values is low, and the global solution is the destination for a high percentage of the initial feasible points, as shown in Figure \ref{fig:nonConvexLess}. We consider this the easier non-convex case.
    \item There are many different terminal values, and the global solution is obtained from a lower percentage of initial points, as shown in Figure \ref{fig:nonConvexMore}. This is the harder non-convex case.
\end{itemize}

~~\\
Interestingly, and perhaps not too surprisingly, the cases in the first category are among those that Gurobi could also solve very quickly. Instances in the second group were solved by Gurobi in good time (1-100 seconds), and the instances in the third group took much longer. 

~~\\ We also conducted some sensitivity analysis. Changing the value of $b$, the relative hardness of the second and third categories {\em reversed} for Gurobi, but it does not have much effect on the first category.
For GAMA,  the solving times for all categories are almost equal and independent of the hardness of $Q$ or value of $b$ and depend on the number of initial feasible solutions (seeds) $l$.
~~\\

We further probed how the value of $b$ affects the Gurobi performance. The problem samples that we mentioned earlier all had $b$ values of either $10$ or $40$, which is equally ($=15$) deviated from the center to the right or the left. From the Gurobi solver point of view, in almost all cases, switching the $b$ value to the opposite side of the center, for easier problems makes the problem harder, and turns harder problems into easier ones (whereas in our solver all of the problems are solved in the same amount of time). Changing the $b$ value from either $10$ or $40$ to the middle value ($25$) transforms the nonconvex problems and even some of the convex problems into extremely difficult; even after several hours, Gurobi was not able to solve any of them optimally (whereas our algorithm still solves them in about $1sec$ on average). We speculate that the reason for this degradation in Gurobi's performance is that as we move the value of $b$ closer to the middle point, the total number of possible nodes in the constraint polytope becomes maximum at $b = n/2$: $max \left( {\begin{array}{*{20}{c}}
  n \\ 
  b 
\end{array}} \right) = \left( {\begin{array}{*{20}{c}}
  n \\ 
  n/2 
\end{array}} \right)$ Consequently, the number of degenerate solutions also increases. This is a hurdle for the Gurobi solver or other such solvers that return only one optimal solution, whereas \textit{our method returns all the degenerate solutions} reached after augmentation terminations.

\subsection{Quadratic Semi-Assignment Problems 1 (QSAP1)}
Here we consider nonconvex nonlinear combinatorial problems with a quadratic term as the objective function and $k$ separate {\em horizontal} cardinality constraints. Quadratic Semi-Assignment Problems (QSAP) have many applications \cite{pitsoulis_quadratic_2009}. The problem is of the form:
\begin{equation}
\begin{array}{*{20}{c}}
  {(P2)}&{\min \left\{ {{\mathbf{c}^T}X + {X^T}QX:\left( {{I_n} \otimes {\mathbf{1}}_k^T} \right)X = {\mathbf{b}}} \right\}} 
\end{array}
\end{equation}
where ${\mathbf{b}} \in \mathbb{Z}_ + ^n$ and
${\mathbf{x}}_i^T = \left[ {\begin{array}{*{20}{c}}
  {{x_{i,1}}}&{{x_{i,2}}}& \cdots &{{x_{i,k}}} 
\end{array}} \right] \in {\mathbb{B}^k}$
is the $k \times 1$ column vector depicting the $i^{th}$ assignment (aka $i^{th}$ brick) and 
\begin{equation}\label{eq:X}
{X^T} = \left[ {\begin{array}{*{20}{c}}
  {{\mathbf{x}}_1^T}&{{\mathbf{x}}_2^T}& \cdots &{{\mathbf{x}}_n^T} 
\end{array}} \right] \in {\mathbb{B}^{kn}}
\end{equation}
is the concatenation of all assignment vectors (bricks).
\subsubsection{Feasible solutions}
The Graver basis of $A = I_n \otimes \mathbf{1}_k^T$ can be generated using equation (\ref{eq:gravermanyone}). To solve the problem using our multi seeded approach, we need to generate many (assume $l$) uniformly distributed feasible solutions of the linear constraint ${\left( {{I_n} \otimes {\mathbf{1}}_k^T} \right)x = {\mathbf{b}}}$, $b_i \leqslant n$. 
The total number of solutions is $\prod\limits_{i = 1}^n {\left( {\begin{array}{*{20}{c}}
  k \\ 
  {{b_i}} 
\end{array}} \right)}$, which can be very large. As before, we want to generate only $l$ uniformly distributed number of them as initial feasible solutions. This is done in two stages:
\begin{itemize}
    \item For each of the $l$ solutions, that has $n$ sections (bricks), we  randomly choose a section among $n$, then generate a size $k$ vector with $1$ placed in  $b_i$ locations with indices \textit{sampled uniformly at random} from $1$ to $k$ (without replacement), and $0$ placed in other locations. 
    \item We uplift this solution to randomly chosen section $i$, in an $nk$ size (initially set to zero) vector. 
    \end{itemize}

\subsubsection{QSAP1 results} \label{subsec:QSAP1results}

\begin{figure}[H]
\centering
\includegraphics[width=16cm]{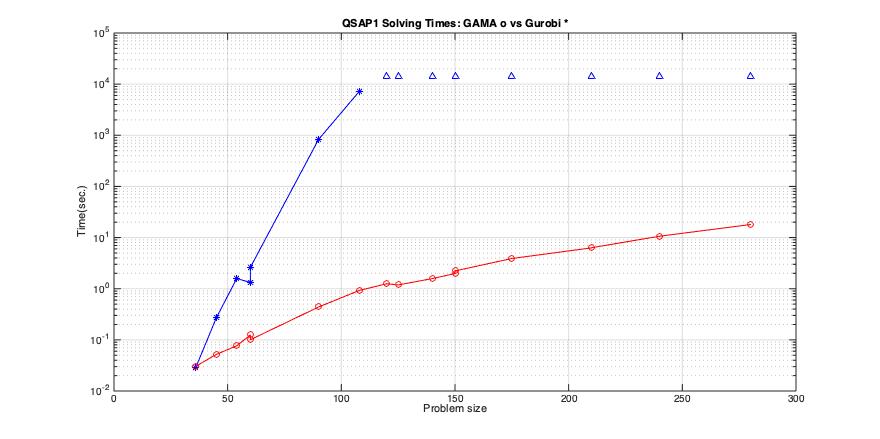}
\caption{Solving time of various size QSAP1 problem instances: GAMA 'o' vs. Gurobi '*'}
\label{fig:timesQSAP1}
\end{figure}

We have three batches of tests.
\begin{itemize}
    \item 
We generated $10$ QSAP problem instances\footnote{We used the QSAP problem instance generator provided by D.E. Bernal.} of size $k \times n = 30 \times 10 = 300$. 
Considering the number of initial feasible points to be $l = 100$, our algorithm solves all problems with an average time of $t_{av} = 7.148sec$. 
\item We increased the size of the $10$ QSAP problems to $k \times n = 30 \times 30 = 900$, and our algorithm solved them with an average time of $t_{av} = 132.5sec = 2.2min$.
The same sets of problems (sizes $300$ and $900$) were passed to the Gurobi MIP solver and none were solved after more than $8$ hours.

\item We used a set of problem instances with various sizes from small to large, from the CMU QSAP problem instance generator. This set consists of $17$ problem instances of  $\{ (k \times n)\}$ sizes:
\begin{equation}\nonumber
\left\{ {\begin{array}{*{20}{c}}
  {\underbrace{3 \times 12}_{36}}&{\underbrace {3 \times 15}_{45}}&{\underbrace {3 \times 18}_{54}}&{\underbrace {3 \times 20}_{60}}&{\underbrace {4 \times 15}_{60}}&{\underbrace {5 \times 18}_{90}}&{\underbrace {6 \times 18}_{108}}&{\underbrace {6 \times 20}_{120}}&{\underbrace {5 \times 25}_{125}}&{\underbrace {4 \times 35}_{140}}&{\underbrace {5 \times 30}_{150}}&{\underbrace {6 \times 25}_{150}}&{\underbrace {7 \times 25}_{175}}&{\underbrace {6 \times 35}_{210}}&{\underbrace {7 \times 30}_{210}}&{\underbrace {8 \times 30}_{240}}&{\underbrace {8 \times 35}_{280}} 
\end{array}} \right\}
\end{equation}
Here we chose the number of starting feasible points equal to the size of the problem $l = k \times n$.
The results are shown in Figure \ref{fig:timesQSAP1}. 
\end{itemize}

~~\\
GAMA optimally solved all $17$ samples in
\begin{equation}\nonumber
\left(\begin{array}{ccccccccccccccccc} 0.03034 & 0.0521 & 0.07652 & 0.1248 & 0.1007 & 0.4394 & 0.929 & 1.251 & 1.186 & 1.578 & 1.983 & 2.227 & 3.898 & 6.288 & 6.329 & 10.63 & 17.97 \end{array}\right)
\end{equation}
seconds, in a total of $55.1sec$, with an average time of $t_{av} = 3.2407sec$. The Gurobi solver, on the other hand, could only finish seven of them in
$\left(\begin{array}{ccccccc} 0.02898 & 0.2686 & 1.585 & 1.323 & 2.629 & 831.3 & 7297.0 \end{array}\right)$
seconds. The other 10 instances could not be solved, and after $4\ hours$ for each instance, the Gurobi solver was terminated. On the largest size that Gurobi could solve, the $7^{th}$ problem instance, GAMA had a $7850+$ factor speedup. 

~~\\

An observation that we can have from Figure \ref{fig:timesQSAP1} is that, for problems of a similar complexity level (roughly like our problem instances in QSAP1), a linear increase in the problem size results in a linear increase in the {\em logarithm} of solution time, which means an exponential increase in time for both GAMA  and Gurobi. As can be seen quite clearly, however, the slope of Gurobi's line is much higher than that of GAMA.
~~\\

\begin{figure}[H]
\centering
\includegraphics[width=16cm]{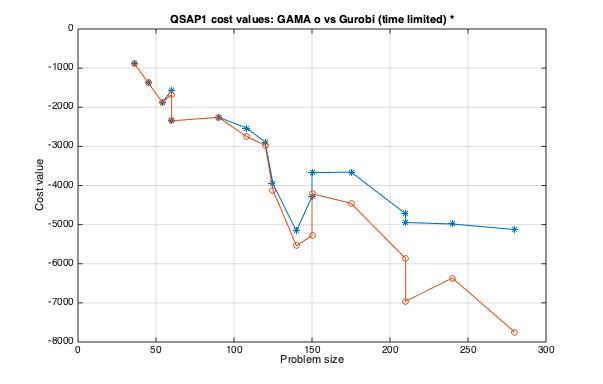}
\caption{Cost values vs. size of time limited Gurobi '*' and GAMA 'o'.}
\label{fig:timelimited}
\end{figure}

It is well known that the {\em total} Gurobi time, like that of many best-in-class exact MIP solvers, consists of separate phases: After finding a feasible solution, they enter the phase of branch-and-bound search and \textit{improving} the incumbent solution, and then switch into the phase of {\em proving} optimality\footnote{Recall that in GAMA, assuming that we have at least one feasible solution in the deepest \textit{convex} region, the termination of Graver augmentation path is a proof of the solution optimality.} \cite{berthold_feasibility_2018}. (The solver may switch back and forth between \textit{improving} and \textit{proving} phases dynamically to minimize the overall time of the solving process.) To make our comparisons on more equal footing, therefore, we also  evaluated the quality of solutions acquired by Gurobi (before spending time on proving optimality). To do this, we set the Gurobi's \textit{timelimit} parameter to three times the average GAMA solving time ($3\times t_{av}$).  We find that for small sized QSAP1 problems the results match the optimal, but for larger sizes the solution obtained within this time limit degrades rapidly (Figure \ref{fig:timelimited}). Consequently, we believe that GAMA finds better solutions faster as the problem size increases.
\subsection{Quadratic Semi-Assignment Problems 2 (QSAP2)}\label{subsec:QSAP2}
Here we solve the nonconvex nonlinear combinatorial problem with a quadratic term as the objective function and $n$ separate {\em vertical} cardinality constraints. Quadratic Semi-Assignment Problems 2 (QSAP2) \cite{pitsoulis_quadratic_2009} is of the form:
\begin{equation}\label{eq:QSAP2}
\begin{array}{*{20}{c}}
  {(P3)}&{\min \left\{ {{\mathbf{c}^T}X + {X^T}QX:\left( {{\mathbf{1}}_n^T \otimes {I_k}} \right)X = {\mathbf{b}}} \right\}} 
\end{array}
\end{equation}
where ${\mathbf{b}} \in \mathbb{Z}_ + ^k$ and similarly 
${\mathbf{x}}_i^T = \left[ {\begin{array}{*{20}{c}}
  {{x_{i,1}}}&{{x_{i,2}}}& \cdots &{{x_{i,k}}} 
\end{array}} \right] \in {\mathbb{B}^k}$
is the $k \times 1$ column vector depicting $i^{th}$ assignment (aka $i^{th}$ brick) and 
${X^T} = \left[ {\begin{array}{*{20}{c}}
  {{\mathbf{x}}_1^T}&{{\mathbf{x}}_2^T}& \cdots &{{\mathbf{x}}_n^T} 
\end{array}} \right] \in {\mathbb{B}^{kn}}$
is concatenation of all assignment vectors (bricks).
\subsubsection{Feasible solutions}
The Graver basis of $A = \mathbf{1}_n^T \otimes I_k$ can be generated using equation (\ref{eq:gravermanyeye}). To solve the problem using our multi seeded approach, we generate many (assume $l$) uniformly distributed feasible solutions of the linear constraint ${\left( {{\mathbf{1}}_n^T} \otimes {I_k} \right)X = {\mathbf{b}}}$, $b_i \leqslant n$. 
The total number of solutions is $\prod\limits_{i = 1}^k {\left( {\begin{array}{*{20}{c}}
  n \\ 
  {{b_i}} 
\end{array}} \right)}$, which can be very large. We generate a much smaller number, $l$, uniformly distributed  initial feasible solutions. As before, this is done in two stages. For each of the $l$ solutions, that has $k$ subsections, we initially randomly choose a section among $k$, then generate a size $n$ vector with $1$ placed in  $b_i$ locations with indices \textit{sampled uniformly at random} from $1$ to $n$ (without replacement), and $0$ placed in other locations. Next, we spread the solutions element $k$ apart by  $k$ elements for each random subsection $i$, and place them in an $nk$ size (initially set to zero) vector. 

\subsubsection{QSAP2 results}
We retained the $Q's$ from QSAP1 instances, and using the same $k$ and $n$, we generated $A$ and $b$ based on QSAP2 formulation, equation (\ref{eq:QSAP2}).
~~\\

\begin{figure}[H]
\centering
\includegraphics[width=16cm]{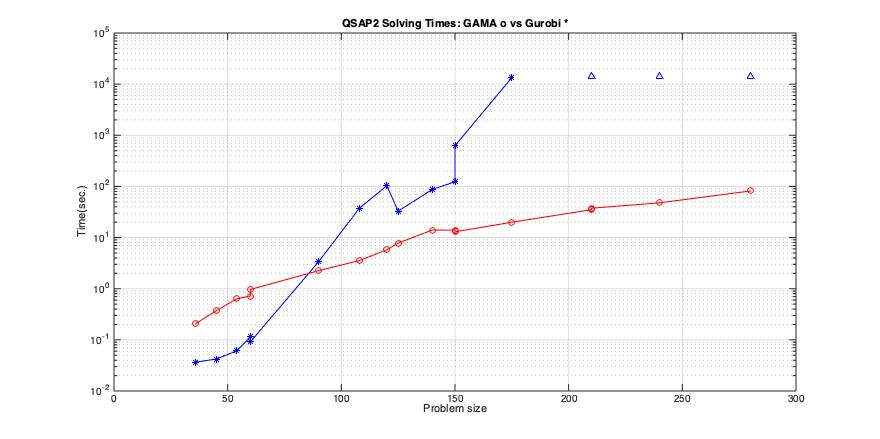}
\caption{Solving time of various size QSAP2 problem instances: GAMA 'o' vs. Gurobi '*'}
\label{fig:timesQSAP2}
\end{figure}

The same set of $17$ problem instances shown in section \ref{subsec:QSAP1results} is used in our testing. GAMA solved all of them with times that ranged from $0.208sec.$ to $81.53sec.$
Gurobi  solved $13$ of the $17$ instances. The other $4$ instances could not be solved in $4\ hours$. 
The results are shown in Figure \ref{fig:timesQSAP2}. The last four instances not completed by Gurobi after 4 hours each, are shown by $\vartriangle$ at the 4 hour border time. 
In the largest instance that Gurobi did solve (in $3.72\ hours$), our method is $650$ times faster. 
As before, the linear dependency between the logarithm of time versus problem size exists here as well, with almost similar slope differences between GAMA and Gurobi. A difference is that, here in QSAP2, the initial crossing point is on higher problem sizes than in QSAP1, indicating that for a much larger range of smaller problems Gurobi is faster, but as the size increases, GAMA outperforms significantly.

\subsection{Quadratic Assignment Problem (QAP)}
Here we solve the nonconvex nonlinear combinatorial problem with a quadratic term as the objective function and $k$ separate horizontal cardinality constraints in addition to $n$ separate vertical cardinality constraints. The Quadratic Assignment Problems (QAP) is of the form:
\begin{equation}
\begin{array}{*{20}{c}}
  {(P4)}&{\min \left\{ {{\mathbf{c}^T}X + {X^T}QX:\left( {\left( {{\mathbf{1}}_n^T \otimes {I_k}} \right) \oplus \left( {{I_n} \otimes {\mathbf{1}}_k^T} \right)} \right)X = {\mathbf{b}}} \right\}}, 
\end{array}
\end{equation}
where ${\mathbf{b}} \in \mathbb{Z}_ + ^{k+n}$ and, similarly, 
${\mathbf{x}}_i^T = \left[ {\begin{array}{*{20}{c}}
  {{x_{i,1}}}&{{x_{i,2}}}& \cdots &{{x_{i,k}}} 
\end{array}} \right] \in {\mathbb{B}^k}$
is the $k \times 1$ column vector depicting the $i^{th}$ assignment (aka $i^{th}$ brick) and 
${X^T} = \left[ {\begin{array}{*{20}{c}}
  {{\mathbf{x}}_1^T}&{{\mathbf{x}}_2^T}& \cdots &{{\mathbf{x}}_n^T} 
\end{array}} \right] \in {\mathbb{B}^{kn}}$
is the concatenation of all assignment vectors (bricks).
\subsubsection{Connection to Random Binary Matrices}
The Graver basis of $A = {\left( {{\mathbf{1}}_k^T \otimes {I_n}} \right) \oplus \left( {{I_k} \otimes {\mathbf{1}}_n^T} \right)}$ can be generated by the procedure described in subsection (\ref{ss:combgraver}).

~~\\
Calculation of randomly generated feasible solutions of $\mathbf{A}X = \mathbf{b}$ constraints for QAP problem categories is not as straightforward as in the previous cases. 
We connect this problem to random binary matrix theory and observe that some known algorithms in that area can be adapted. Additionally, we propose a novel approach -- also based on the Graver basis (Null basis can also be used) -- that not only helps us find random feasible solutions for QAP, but can also be applied to study random binary matrices as well.

~~\\
We rearrange the main vector $X$ defined in equation (\ref{eq:X}), such that $n$ size $k$ sub-vectors $\mathbf{x}_i$ become columns of a $k \times n$ matrix $X_M$:
\begin{equation} \label{eq:XM}
{X_M} = \left[ {\begin{array}{*{20}{c}}
  {{{\mathbf{x}}_1}}&{{{\mathbf{x}}_2}}& \cdots &{{{\mathbf{x}}_n}} 
\end{array}} \right] \in {\mathbb{B}^{k \times n}}
\end{equation}
where $X = vec\left( {{X_M}} \right)$, and $vec$ is the vectorization\footnote{$vec\left( {{X_M}} \right) = \sum\limits_{i = 1}^n {{e_i} \otimes X} {e_i}$} operator.
Thus, the QAP constraint
$$\left( {\left( {{\mathbf{1}}_n^T \otimes {I_k}} \right) \oplus \left( {{I_n} \otimes {\mathbf{1}}_k^T} \right)} \right)X = {\mathbf{b}}$$
becomes
\[\left\{ \begin{gathered}
  {X_M}{{\mathbf{1}}_n} = r \hfill \\
  {\mathbf{1}}_k^T{X_M} = c^T \hfill \\ 
\end{gathered}  \right.\begin{array}{*{20}{c}}
  ,&{\left[ {\begin{array}{*{20}{c}}
  r \\ 
  c 
\end{array}} \right] = {\mathbf{b}},} 
\end{array}\]
which states that row sum ($r$) and column sum ($c$) of matrix $X_M$ should be fixed. This then helps create ${\mathbf{b}}$.

~~\\

The problem of finding random feasible solutions of QAP problems becomes the problem of generating  random binary matrices with fixed row sum and column sum, also known as \textit{doubly stochastic binary matrices}. This problem is addressed in the literature\footnote{Sampling of zero one matrices has many applications in statistical analysis of many fields of study including co-occurrence matrices in evolutionary studies and ecology, multivariate binary time series, affiliation matrices in sociology, and item response analysis in psycho-metrics.}, and  several known algorithms tackle it.  Rycer \cite{rycer_matrices_1960} was the first to study such matrices. Some combinatorial properties of Rycer matrices were studied in \cite{brualdi_matrices_1980} by connecting them to bipartite matrices and hyper-graphs. Enumeration of binary matrices with known row and column sum is also studied \cite{wang_precise_1998}.
The {\em Curveball} algorithm (or its reincarnation \cite{strona_fast_2014}) is known to be a fast method that creates uniform random samples \cite{carstens_proof_2015} of binary matrices with fixed row and column sums. Other approaches for exact counting and sampling of binary matrices with specified sums based on dynamic programming have been devised \cite{miller_exact_2013}.

\subsubsection{Feasible solutions: A novel Graver basis approach}
Here we describe a novel approach that we have used to find a random distribution of binary matrices, coincidentally based on the Graver basis. Similar to the swap in the Curveball algorithm (and some other algorithms), we start with one single solution and add a randomly chosen Graver element to it (while checking for lower bound $0$ and upper bound $1$ on all terms) and repeat this a randomly chosen number of times, to reach the next feasible solution. Repeating this procedure, again and again, generates more feasible solutions. 
That is, let $X_0$ be the initial feasible solution and $Card\left( {\mathcal{G}({\mathbf{A}})} \right) = N$. We select from $1...N$ a random number ($n_r$), then chose $n_r$ random indices from $1...N$, add them to $X_0$ to reach to the $X_1$, and repeat.
\subsubsection{QAP results}
  The constraint matrix $\textbf{A}$ is generated based on the QAP. Values for $\mathbf{b}$ are generated by initially generating a random $k \times n$ binary matrix and using its row sum and column sum accordingly, to guarantee feasibility such that the sum of column vectors equals the sum of row vectors. The same set of $17$ problem instances described in section \ref{subsec:QSAP1results} is used.

~~\\

\begin{figure}[H]
\centering
\includegraphics[width=15cm]{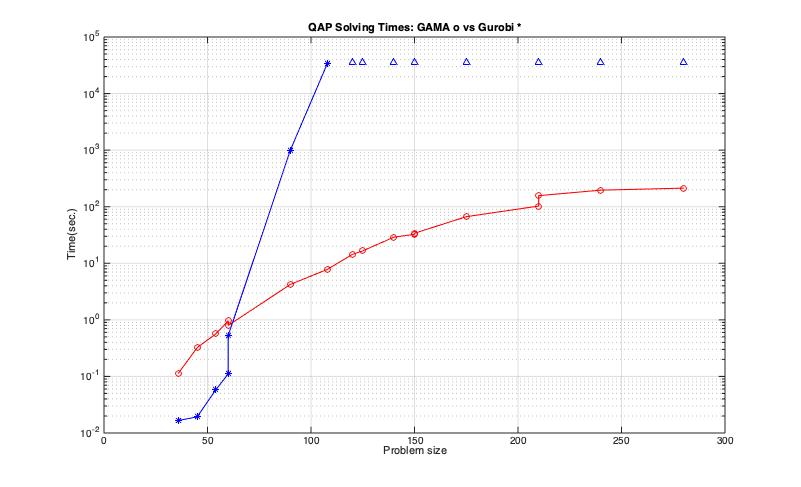}
\caption{Solving time of various size QAP problem instances: GAMA 'o' vs. Gurobi '*'}
\label{fig:timesQAP}
\end{figure}

GAMA solved all 17 instances with times ranging from $0.1144sec.$ to $212.86sec$.
Gurobi solved $7$ of them, and could not solve the other 10 even after $10\ hours$ for each. Results are shown in Figure \ref{fig:timesQAP}. The instances not completed by Gurobi after $10$ hours are shown by $\vartriangle$ at the $10$ hour border time. 
In the largest instance that Gurobi solved in under
$10$ hours (problem size $108$ in $9.6\ hours$), GAMA ($7.8sec.$) is $4407$ times faster. 

~~\\
As can be seen in Figure \ref{fig:timesQAP}, the same linear dependency between logarithm of time versus problem size exists here. Gurobi could solve small sized QAP problem instances faster than GAMA, but the  slope is much higher than before, due to tighter double (row and column) constraints.
~~\\

Further speedup of GAMA is possible. A study of Graver bases elements that are used in the augmentation paths indicates that the majority of bases that cause cost reduction are from basis elements that are the result of lower brick number liftings ($t=2\sim4$). In other words, the higher $t$-cycle equivalents of the Hilbert basis are not used frequently. In cases when $k$ and $n$ are both large, thus the number of Graver basis elements is large, instead of generating all the liftings of a Graver set, we generate only liftings with a lower number of bricks and created a fixed Graver basis repository, randomly generated elements from the higher size liftings during augmentation, and used them on the fly. This procedure and similar modification of it limits the memory requirements substantially. It is important to note that this selective random generation is possible only when we systematically create the Graver basis like we do, which is not the approach using a completion procedure such as the Pottier algorithm \cite{pottier_euclidean_1996}.

\section{Conclusions}
In this paper, we proposed and tested a novel method, the Graver Augmented Multi-seeded Algorithm (GAMA), for  non-convex, non-linear integer programs. For problem classes, that include Cardinality Boolean Quadratic Programs (CBQP), Quadratic Semi-Assignment Programs, and Quadratic Assignment Programs, we develop procedures for (1) systematically constructing Graver basis elements and (2) finding many feasible solutions that are uniformly spread out. We performed extensive numerical testing on instances that arise in industry, and conducted a sensitivity analysis to understand why GAMA vastly outperforms existing best-in-class commercial solvers.   

\section*{Acknowledgments}
\addcontentsline{toc}{section}{Acknowledgements}
The authors thank D. E. Bernal and Professor I. Grossmann from CMU for providing us with instances of CBQP and QSAP, and for providing helpful feedback on an earlier draft. 

\section*{References}
\addcontentsline{toc}{section}{References}
\label{sec:bibliography}
\bibliographystyle{amsalpha}
\bibliography{main.bib}

\newcommand{\etalchar}[1]{$^{#1}$}
\providecommand{\bysame}{\leavevmode\hbox to3em{\hrulefill}\thinspace}
\providecommand{\MR}{\relax\ifhmode\unskip\space\fi MR }
\providecommand{\MRhref}[2]{%
  \href{http://www.ams.org/mathscinet-getitem?mr=#1}{#2}
}
\providecommand{\href}[2]{#2}
\begin{thebibliography}{DLHOW08}

\bibitem[ADT19]{alghassi_graver_2019}
Hedayat Alghassi, Raouf Dridi, and Sridhar Tayur, \emph{Graver {Bases} via
  {Quantum} {Annealing} with {Application} to {Non}-{Linear} {Integer}
  {Programs}}, arXiv:1902.04215 [quant-ph] (2019), arXiv: 1902.04215.

\bibitem[BEHM06]{bruglieri_annotated_2006}
Maurizio Bruglieri, Matthias Ehrgott, Horst~W. Hamacher, and Francesco
  Maffioli, \emph{An annotated bibliography of combinatorial optimization
  problems with fixed cardinality constraints}, Discrete Applied Mathematics
  \textbf{154} (2006), no.~9, 1344--1357.

\bibitem[BHK18]{berthold_feasibility_2018}
Timo Berthold, Gregor Hendel, and Thorsten Koch, \emph{From feasibility to
  improvement to proof: three phases of solving mixed-integer programs},
  Optimization Methods and Software \textbf{33} (2018), no.~3, 499--517.

\bibitem[Bil05]{billionnet_different_2005}
Alain Billionnet, \emph{Different {Formulations} for {Solving} the {Heaviest}
  {K}-{Subgraph} {Problem}}, INFOR: Information Systems and Operational
  Research \textbf{43} (2005), no.~3, 171--186.

\bibitem[BLSR99]{bigatti_computing_1999}
A.~M. Bigatti, R.~La~Scala, and L.~Robbiano, \emph{Computing {Toric} {Ideals}},
  Journal of Symbolic Computation \textbf{27} (1999), no.~4, 351--365.

\bibitem[BPT00]{bertsimas_new_2000}
Dimitris Bertsimas, Georgia Perakis, and Sridhar Tayur, \emph{A {New}
  {Algebraic} {Geometry} {Algorithm} for {Integer} {Programming}}, Management
  Science \textbf{46} (2000), no.~7, 999--1008.

\bibitem[Bru80]{brualdi_matrices_1980}
Richard~A. Brualdi, \emph{Matrices of zeros and ones with fixed row and column
  sum vectors}, Linear Algebra and its Applications \textbf{33} (1980),
  159--231.

\bibitem[Car15]{carstens_proof_2015}
C.~J. Carstens, \emph{Proof of uniform sampling of binary matrices with fixed
  row sums and column sums for the fast {Curveball} algorithm}, Physical Review
  E \textbf{91} (2015), no.~4, 042812.

\bibitem[CT91]{conti_buchberger_1991}
Pasqualina Conti and Carlo Traverso, \emph{Buchberger {Algorithm} and {Integer}
  {Programming}}, Proceedings of the 9th {International} {Symposium}, on
  {Applied} {Algebra}, {Algebraic} {Algorithms} and {Error}-{Correcting}
  {Codes} (London, UK, UK), {AAECC}-9, Springer-Verlag, 1991, pp.~130--139.

\bibitem[DLHO{\etalchar{+}}09]{de_loera_convex_2009}
J.~A. De~Loera, R.~Hemmecke, S.~Onn, U.~G. Rothblum, and R.~Weismantel,
  \emph{Convex integer maximization via {Graver} bases}, Journal of Pure and
  Applied Algebra \textbf{213} (2009), no.~8, 1569--1577.

\bibitem[DLHOW08]{de_loera_n-fold_2008}
Jesús~A. De~Loera, Raymond Hemmecke, Shmuel Onn, and Robert Weismantel,
  \emph{N-fold integer programming}, Discrete Optimization \textbf{5} (2008),
  no.~2, 231--241.

\bibitem[GO19]{gurobi}
LLC Gurobi~Optimization, \emph{Gurobi optimizer reference manual}, 2019.

\bibitem[Gra75]{graver_foundations_1975}
Jack~E. Graver, \emph{On the foundations of linear and integer linear
  programming {I}}, Mathematical Programming \textbf{9} (1975), no.~1, 207--226
  (en).

\bibitem[HOW11]{hemmecke_polynomial_2011}
Raymond Hemmecke, Shmuel Onn, and Robert Weismantel, \emph{A polynomial
  oracle-time algorithm for convex integer minimization}, Mathematical
  Programming \textbf{126} (2011), no.~1, 97--117 (en).

\bibitem[HS95]{hosten_grin:_1995}
Serkan Hoşten and Bernd Sturmfels, \emph{{GRIN}: {An} implementation of
  {Gröbner} bases for integer programming}, Integer {Programming} and
  {Combinatorial} {Optimization}, Lecture {Notes} in {Computer} {Science},
  Springer, Berlin, Heidelberg, May 1995, pp.~267--276 (en).

\bibitem[HS07]{hosten_finiteness_2007}
Serkan Hoşten and Seth Sullivant, \emph{A finiteness theorem for {Markov}
  bases of hierarchical models}, Journal of Combinatorial Theory, Series A
  \textbf{114} (2007), no.~2, 311--321.

\bibitem[LG17]{lima_solution_2017}
Ricardo~M. Lima and Ignacio~E. Grossmann, \emph{On the solution of nonconvex
  cardinality {Boolean} quadratic programming problems: a computational study},
  Computational Optimization and Applications \textbf{66} (2017), no.~1, 1--37
  (en).

\bibitem[LORW12]{lee_quadratic_2012}
Jon Lee, Shmuel Onn, Lyubov Romanchuk, and Robert Weismantel, \emph{The
  quadratic {Graver} cone, quadratic integer minimization, and extensions},
  Mathematical Programming \textbf{136} (2012), no.~2, 301--323 (en).

\bibitem[MH13]{miller_exact_2013}
Jeffrey~W. Miller and Matthew~T. Harrison, \emph{{EXACT} {SAMPLING} {AND}
  {COUNTING} {FOR} {FIXED}-{MARGIN} {MATRICES}}, The Annals of Statistics
  \textbf{41} (2013), no.~3, 1569--1592.

\bibitem[MSW04]{murota_optimality_2004}
Kazuo Murota, Hiroo Saito, and Robert Weismantel, \emph{Optimality criterion
  for a class of nonlinear integer programs}, Operations Research Letters
  \textbf{32} (2004), no.~5, 468--472.

\bibitem[Onn10]{onn_nonlinear_2010}
Shmuel Onn, \emph{Nonlinear {Discrete} {Optimization}: {An} {Algorithmic}
  {Theory}}, European Mathematical Soc., 2010 (en), Google-Books-ID:
  kcD5DAEACAAJ.

\bibitem[Pit09]{pitsoulis_quadratic_2009}
Leonidas Pitsoulis, \emph{Quadratic {Semi}-assignment {Problem}}, Encyclopedia
  of {Optimization} (Christodoulos~A. Floudas and Panos~M. Pardalos, eds.),
  Springer US, Boston, MA, 2009, pp.~3170--3171 (en).

\bibitem[Pot96]{pottier_euclidean_1996}
Loïc Pottier, \emph{The {Euclidean} {Algorithm} in {Dimension} {N}},
  Proceedings of the 1996 {International} {Symposium} on {Symbolic} and
  {Algebraic} {Computation} (New York, NY, USA), {ISSAC} '96, ACM, 1996,
  pp.~40--42.

\bibitem[Ryc60]{rycer_matrices_1960}
Herbert~John Rycer, \emph{Matrices of zeros and ones.}, Bulletin of the
  American Mathematical Society \textbf{66} (1960), no.~6, 442--464.

\bibitem[SNB{\etalchar{+}}14]{strona_fast_2014}
Giovanni Strona, Domenico Nappo, Francesco Boccacci, Simone Fattorini, and
  Jesus San-Miguel-Ayanz, \emph{A fast and unbiased procedure to randomize
  ecological binary matrices with fixed row and column totals}, Nature
  Communications \textbf{5} (2014), 4114 (en).

\bibitem[SS03]{santos_higher_2003}
Francisco Santos and Bernd Sturmfels, \emph{Higher {Lawrence} configurations},
  Journal of Combinatorial Theory, Series A \textbf{103} (2003), no.~1,
  151--164.

\bibitem[ST97]{sturmfels_variation_1997}
Bernd Sturmfels and Rekha~R. Thomas, \emph{Variation of cost functions in
  integer programming}, Mathematical Programming \textbf{77} (1997), no.~2,
  357--387 (en).

\bibitem[TTN95]{tayur_algebraic_1995}
Sridhar~R. Tayur, Rekha~R. Thomas, and N.~R. Natraj, \emph{An algebraic
  geometry algorithm for scheduling in presence of setups and correlated
  demands}, Mathematical Programming \textbf{69} (1995), no.~1-3, 369--401
  (en).

\bibitem[WZ98]{wang_precise_1998}
Bo-Ying Wang and Fuzhen Zhang, \emph{On the precise number of (0,1)-matrices in
  {A}({R},{S})}, Discrete Mathematics \textbf{187} (1998), no.~1, 211--220.

\end{thebibliography}

\end{document}